\documentclass[11pt]{amsart}
\usepackage{amssymb}
\usepackage{amsmath}
\usepackage{amscd}
\usepackage[dvips]{color}

\newtheorem{prop-def}{Proposition-Definition}[section]

\begin{document}

\title{Symplectic structures on $3$-Lie algebras}

\author{RuiPu  Bai}
\address{College of Mathematics and Computer Science,
Hebei University, Baoding 071002, P.R. China}
\email{bairuipu@hbu.edu.cn}

\author{ Shuangshuang Chen}
\address{College of Mathematics and Computer Science,
Hebei University, Baoding 071002, P.R. China}
\email{chenss0416@13.com}

\author{Rong Cheng}
\address{College of Mathematics and Computer Science,
Hebei University, Baoding 071002, P.R. China}
\email{chengrongbaoding@126.com}

\date{}
\maketitle

\thispagestyle{plain} \markright{\upshape \hfill \sc Symplectic
structures on $3$-Lie algebras \hfill}

\vspace{-5mm}
\begin{abstract} The symplectic structures on $3$-Lie algebras and metric symplectic $3$-Lie algebras are studied. For arbitrary $3$-Lie algebra $L$,
 infinite many metric symplectic $3$-Lie algebras are constructed.
It is proved that a metric $3$-Lie algebra $(A, B)$ is a metric symplectic $3$-Lie algebra if and only if there exists an invertible
derivation $D$ such that $D\in Der_B(A)$,  and is also proved that every   metric symplectic $3$-Lie algebra $(\tilde{A}, \tilde{B}, \tilde{\omega})$ is a $T^*_{\theta}$-extension of a   metric symplectic $3$-Lie algebra  $(A, B, \omega)$. Finally, we construct a   metric symplectic double extension of a  metric symplectic $3$-Lie algebra by means of a special derivation.

\end{abstract}

\baselineskip=18pt

\section{Introduction}
\setcounter{equation}{0}
\renewcommand{\theequation}
{1.\arabic{equation}}

The notion of $3$-Lie algebra was introduced in \cite{F}. It is a vector space with a ternary linear skew-symmetric  multiplication satisfying the generalized Jacobi identity (or Filippov identity). $3$-Lie algebras, especially, metric $3$-Lie algebras are applied in many fields in mathematics and mathematical physics.
 Motivated by some problems of quark dynamics, Nambu \cite{N} introduced a $3$-ary generalization of Hamiltonian dynamics by means of the $3$-ary
Poisson bracket
$$
[f_1, f_2, f_3]= \det\Big(\frac{\partial f_i}{\partial x_j}\Big)
$$
which satisfies the generalized Jacobi identity
$$
[[f_1, f_2, f_3], g_2, g_2]= [[f_1, g_2, g_3], f_2, f_3]+[f_1, [f_2, g_2, g_3], f_3]+[f_1, f_2, [f_3, g_2, g_3]].
$$
Following this line, Takhtajan   \cite{T} developed systematically
 the foundation of the theory of $n$-Poisson or
Nambu-Poisson manifolds.
Metric $3$-Lie algebras are applied to the study of the
supersymmetry and gauge symmetry transformations of the world-volume
theory of multiple coincident M2-branes; the Bagger-Lambert theory
has a novel local gauge symmetry which is based on a metric $3$-Lie
algebra \cite{BL, HHM}. The generalized Jacobi identity can be
regarded as a generalized Plucker relation in the physics
literature \cite{HCK,G,P}.

Authors in \cite{bai3} studied the structure of metric $n$-Lie algebras. It is  an $n$-Lie algebra with a non-degenerate $ad$-invariant symmetric bilinear form.
The ordinary gauge theory requires a positive-definite metric to guarantee that the theory
possesses positive-definite kinetic terms and to prevent violations of unitarity due to propagating
ghost-like degrees of freedom. But very few metric $n$-Lie algebras admit positive-definite metrics
(see \cite{P,JLZ}); Ho, et al. in \cite{HHM} confirmed that there are no non-strong semisimple $n$-Lie algebras \cite{bai6}
with positive-definite metrics for $n=5,6,7,8$. They also gave examples of 3-Lie algebras whose
metrics are not positive-definite and observed that generators of zero norm are common in $3$-Lie
algebras. Papers \cite{bai4, bai5} studied the module-extension of $3$-Lie algebras and $T_{\theta}^*$-extension of $n$-Lie algebras. So we can obtain more metric $3$-Lie algebras
by $3$-Lie algebras and their modules.

We know that Lie groups which admit a bi-invariant pseudo-Riemannian metric and also a left-invariant
symplectic form are nilpotent Lie groups and their geometry (and, consequently, that
of their associated homogeneous spaces) is very rich. In particular, they carry two left-invariant affine structures: one defined by the symplectic form (which is well-known) and
another which is compatible with a left-invariant pseudo-Riemannian metric. The paper \cite{Bor} studied quadratic Lie algebras over a field $K$ of null characteristic which admit,
at the same time, a symplectic structure. It is proved that if $K$ is algebraically closed every
such Lie algebra may be constructed as the $T^*$
-extension of a nilpotent algebra admitting an invertible derivation and also as the double extension of another quadratic
symplectic Lie algebra by the one-dimensional Lie algebra. In this paper we study the  metric $3$-Lie algebra which, at same time, admits a  symplectic structure. We call it a metric
 symplectic $3$-Lie algebra.

Throughout this paper, $F$ denotes
an algebraically closed field $F$ of characteristic zero. Any
bracket that is not listed in the multiplication of a $3$-Lie
algebra is assumed to be zero.
The symbol $\oplus$ will be frequently used.Unless
other thing is stated, it will only denote the direct sum of vector
spaces.

\section{Fundamental notions}
\setcounter{equation}{0}
\renewcommand{\theequation}
{2.\arabic{equation}}

 {\it A $3$-Lie algebra} \cite{F} is a vector space $L$ over a field F on which  a linear
multiplication  $[~ , ~ , ~ ]: L\wedge L\wedge L \rightarrow L$
satisfying generalized Jacobi identity (or Filippov identity)

$$ [[x_1, x_2, x_3],y_2, y_3]=\sum\limits_{i=1}^3[x_1, \cdots, [x_i,y_2, y_3],
\cdots, x_3], ~ \forall x_1, x_2, x_3, y_2, y_3\in L.$$

\vspace{2mm} A subspace $A$ of $L$ is called a {\it subalgebra} ({\it an ideal }) of $L$ if
$[A, A, A]\subseteq A$ $([A, L, L]\subseteq A )$. If $[A, A, A]=0$ $( [A, A, L]=0 )$, than $A$ is
called {\it an abelian subalgebra} ({\it an abelian ideal}) of $L$.

In particular, the subalgebra
generated by the vectors $[x_1, x_2, x_3]$ for all $x_1,
x_2, x_3\in L$ is called the {\it derived algebra} of $L$,
which is denoted by $L^1$. If $L^1=0$, $L$ is called {\it an
abelian algebra.}

A derivation of a $3$-Lie algebra $L$ is a linear mapping $D: L\rightarrow L$ satisfying
$$D[x, y, z]=[Dx, y, z]+[x, Dy, z]+[x, y, Dz], ~ \forall x, y, z\in L.$$
All the derivations of $L$ is a linear Lie algebra, is denoted by $Der(L)$.

A  $3$-Lie algebra $L$ is said to be {\it simple} if $L^1\neq 0$
and it has no ideals distinct from $0$ and itself.

An ideal $I$ of an $3$-Lie algebra $L$ is called  {\it  nilpotent} \cite{K1},
if $I^{s}=0$ for some $s\geq 0$, where $I^{0}=I$ and $I^{s}$ is
defined as
\begin{equation*}
I^{s}=[I^{s-1},  I, L], ~ \mbox{ for } ~ s\geq 1.
\end{equation*}
In the case $I=L$, $L$ is called a nilpotent $3$-Lie algebra.
The abelian ideal $$Z(L)=\{ x\in L~ |~ [x, L, L]=0~ \}$$ is called {\it the center } of $L$.

 Let
$L$ be a $3$-Lie algebra, $V$ be a vector space, $\rho: L\wedge
L\rightarrow End(V)$ be a linear mapping. The pair $(V,\rho)$ is
called {\it a representation} \cite{K1} (or $V$ is an $L$-module)
of $L$ in $V$  if $\rho$
 satisfies $\forall a_1, a_2, a_3, b_1, b_2\in L,$
 $$
[\rho(a_1, a_2), \rho(b_1, b_2)]=\rho([a_1, a_2, b_1],
b_2)+\rho(b_1, [a_1, a_2, b_2]),
$$
$$
\rho([a_{1}, a_2, a_{3}], b_{1})= \rho(a_{2}, a_{3})\rho(a_1,
b_1)-\rho(a_{1}, a_{3})\rho(a_2, b_1)+\rho(a_{1}, a_{2})\rho(a_3,
b_1).
$$
Then  $(V, \rho)$ is a representation of the
$3$-Lie algebra $L$ if and only if the vector space  $Q=L\oplus V$
is a $3$-Lie algebra in the following multiplication
$$
[a_1+v_1, a_2+v_2, a_3+v_3]=[a_1, a_2, a_3]_L+\rho(a_1,
a_2)(v_3)-\rho(a_1, a_3)(v_2)+\rho(a_2, a_3)(v_1).$$ Therefore,  $A$ is a subalgebra and $V$ is an
abelian ideal of the $3$-Lie algebra $L\oplus V$, respectively.

If $(V, \rho)$ is a representation of the $3$-Lie algebra $L$,
$V^*$ is the dual space of $V$. Then $(V^*, \rho^*)$ is also a
representation of $L$, which is called the dual representation of
$(V, \rho)$,  where
$$\rho^*: L\wedge L\rightarrow End(V^*), ~ \rho^*(a, b)f(c)=-f(\rho(a, b)c), ~\forall  a, b, c\in L, f\in V^*.$$

 For $3$-Lie algebra $L$, the joint representation $(L, ad)$ is
 $$ad: L\wedge
L\rightarrow End(L),~~ ad(x, y)(z)=[x, y, z], \forall x, y, z\in
L.$$ Then we obtain the dual
representation  $ad^*: L\wedge L\rightarrow End(L^*)$,
$$(ad^*(x, y)f)(z)=-f(ad(x,
y)z)=-f([x, y, z]), \forall x, y, z\in L, ~~ f\in L^*.$$

  Let $L$ be  a $3$-Lie algebra, $B: L\times L\rightarrow F$ be a non-degenerate  symmetric bilinear form on $L$. If $B$ satisfies
  $$B([x_{1},x_{2},x_{3}],x_{4})+B(x_{3},
 [x_{1},x_{2},x_{4}])=0, \forall x_{1},x_{2},x_{3},x_{4}\in L, \eqno(2.1)$$ then $B$ is called {\it a metric} on  $3$-Lie algebra
 $L$, and $(L, B)$ is called
 {\it a metric $3$-Lie algebra} \cite{bai3}.

 Let $(L, B)$ be a metric $3$-Lie algebra. Denotes
 $$Der_B(L)=\{ D\in Der(L)~|~B(Dx, y)+B(x, Dy)=0, ~ \forall x,y\in  L \}=Der(L)\cap so(L, B).\eqno(2.2)$$

 Let $W$ be a subspace of a metric $3$-Lie algebra $(L, B)$. {\it
The orthogonal complement} of $W$ is defined by
$$
    W^{\bot}=\{x \in L \mid B(w, x)=0~ \mbox{for all} ~ w\in W\}.
$$
Then $W$ is an ideal if and only if $W^{\bot}$ is  an ideal and $(
W^{\bot})^{\bot}=W$. Notice that $W$ is a minimal ideal if and only
if $W^{\bot}$ is maximal.  If $W\subseteq W^{\bot}$, then $W$ is called {\it isotropic}.

 The subspace $W$ is called {\it nondegenerate} if
$B|_{W\times W}$ is nondegenerate, this is equivalent to  $W\cap
W^{\bot}=0$ or $ L =W\oplus W^{\bot}$.  If an ideal $I$ satisfies $I=I^{\bot}$, then  $I$ is called {\it a completely isotropic ideal}.

 If $L$ does not contain nontrivial
nondegenerate ideals, then $L$ is called {\it
$B$-irreducible.} For a metric $3$-Lie algebra $(L, B)$,
it is not difficult to see $$L^1=[L, L, L]=Z(L)^{\bot}.$$

\section{ Symplectic $3$-Lie algebras}
\setcounter{equation}{0}
\renewcommand{\theequation}
{3.\arabic{equation}}

\vspace{2mm}{\bf Definition 3.1} {\it Let  $L$ be a $3$-Lie
algebra over a field $F$, linear mapping $\omega: L\wedge L \rightarrow
F$ be non-degenerate. If $\omega$ satisfies
$$\sum_{i=1}^4
\omega([x_{1}, \cdots, \hat{x_{i}}, \cdots,
x_{4}],(-1)^{i-1}x_{i})=0, ~\forall x_{i}\in L, i=1, 2, 3,
4,\eqno(3.1)$$ then $\omega$ is called a symplectic structure on
$L$, and $(L,\omega)$ is called a symplectic $3$-Lie algebra.}

 An ideal $I$ of a symplectic $3$-Lie algebra $(L, \omega)$ is called  {\it an lagrangian ideal } if and
only if it coincides with its orthogonal with respect to the form
$\omega$.

\vspace{2mm}  If  there exists a metric $B$ and a symplectic
structure $\omega$ on $3$-Lie algebra $L$, respectively, then $(L,
B,\omega)$ is called {\it a metric symplectic 3-Lie algebra.}

By the above definition, if $(L,\omega)$ is  a symplectic $3$-Lie
algebra, then the dimension of $L$ is even.

\vspace{3mm}{\bf Theorem 3.1} {\it  Let $(L, B)$ be a metric
$3$-Lie algebra. Then there exists a symplectic structure on $L$
if and only if there exists a skew-symmetric invertible derivation
$D\in Der_B(L)$.}

\vspace{2mm}{\bf Proof.}  Let $(L, B, \omega)$ be a symplectic
$3$-Lie algebra. Defines $D: L\rightarrow L$ by
$$B(Dx, y)=\omega(x, y), ~~ \forall x,y\in L. \eqno(3.2)$$
Then $D$ is invertible, and from Eq.(3.1), for $\forall x_1, x_2, x_3,
x_4\in L$,

\vspace{2mm}
   $B([Dx_{1},x_{2},x_{3}],x_{4})+B([x_{1},Dx_{2},x_{3}],x_{4})+B([x_{1},x_{2},Dx_{3}],x_{4})-B(D[x_{1},x_{2},x_{3}],x_{4})$

\vspace{2mm} \noindent$
  =-B([x_{2},x_{3},x_{4}],Dx_{1})+B([x_{1},x_{3},x_{4}],D x_{2})-B([x_{1},x_{2},x_{4}],D x_{3})+B([x_{1},x_{2},x_{3}],Dx_{4})
$

\vspace{2mm}\noindent$=\sum_{i=1}^4 \omega([x_{1}, \cdots, \hat{x_{i}},
\cdots, x_{4}],(-1)^{i-1}x_{i})=0.$

Therefore,  $D$ is a skew-symmetric invertible derivation of $(L,
B)$, that is, $D\in Der_B(L).$

Conversely, if $D\in Der_B(L)$ is invertible. Defines $\omega: L\times L\rightarrow F$ by Eq.(3.2).
Then by the above discussion, $\omega$ is non-degenerate, and
satisfies Eq.(3.1). The result follows. $\Box$

{\bf Remark 1}  One might thus think that every symplectic $3$-Lie algebra $(A,
\omega)$ admitting an invertible derivation which is
skew-symmetric for $\omega$ carries a metric structure; but this
is not the case. Let $A$ be a $4$-dimensional $3$-Lie algebra, the
multiplication in a basis $\{ x_{1},x_{2},x_{3},x_{4}\}$  be
defined by
$$[x_{1},x_{2},x_{4}]=x_{3}.$$

Then  the non-degenerate skew-symmetric bilinear form on $A$ given
by $$\omega(x_{1},x_{4})=\omega(x_{2},x_{3})=1$$ provides a
symplectic structure on $A$, and the linear endomorphism of $A$
given by
$$D(x_{1})=2x_{1},D(x_{2})=-x_{2},D(x_{3})=-x_{3},D(x_{4})=-2x_{4}$$
is a skew-symmetric derivation of $(A, \omega)$. But for every
symmetric bilinear form $B: A\times A\rightarrow F$ satisfying
Eq.(2.1), $B$ satisfies $$B(x_3, x_3)=B(x_3, x_1)=B(x_3,
x_2)=B(x_3, x_4)=0.$$ Therefore, $B$ is degenerated. It follows
that there does not exist metric structure on the $3$-Lie algebra
$A.$

  Under the assumptions of Theorem 3.1,
the skew-symmetric derivation $D\in Der_B(L)$ is also
skew-symmetric with respect to the symplectic form $\omega$ since
for all $x,y\in L$,
$$
\omega(Dx,y)=B(D^{2}x,y)=-B(Dx,Dy)=-\omega(x,Dy).
$$

Now for arbitrary $3$-Lie algebra $L$ and a positive integer $n ( n>2)$, we construct a metric symplectic $3$-Lie
algebra. Let $N$ be the set of all non-negative integers, $$F[t]=\{ f(t)=\sum\limits_{i=0}^ma_it^i~|~ a_i\in F, m\in N\}$$ be the algebra of polynomials over $F$.
We consider
 the tensor product of vector spaces $$L_{n}=L\otimes (tF[t]/t^nF[t]),\eqno(3.3)$$
where $tF[t]/t^nF[t]$ is the quotient space of $tF[t]$ module $t^nF[t]$. Then $L_n$ is a nilpotent 3-Lie algebra in the following multiplication
 $$[x\otimes t^{\bar{p}}, y\otimes t^{\bar{q}}, z\otimes t^{\bar{r}}]=[x, y, z]_{L}\otimes t^{\overline{p+q+r}},
 x, y, z\in L; p, q, r\in N \setminus \{0\}. \eqno(3.4)$$

  Defines  endomorphism $D$ of $L_n$ by
  $$D(x\otimes t^{\overline{p}})=p(x\otimes t^{\overline{p}}), ~  \forall x\in L, ~ p=1, \cdots, n-1.$$  Then $D$ is an
  invertible derivation of the $3$-Lie algebra $L_n$.

Let  $\tilde{L_n}=L_n\oplus L_n^{\ast}$, where $L_n^*$ is the dual space of $L_n$. Then $(\tilde{L_n}, B)$ is a metric $3$-Lie algebra with the multiplication
$$
  [x+f, y+g, z+h]=[x, y, z]_{L_n}+ad^*(y,z)f-ad^*(x,z)g+
  ad^*(x,y)h,\eqno(3.5)
$$
for $x, y, z\in L_n, f, g, h\in L_n^{\ast}$, and the bilinear form
$$B(x+f, y+g)=f(y)+g(x).\eqno(3.6)$$

Defines linear mapping $\tilde{D}: \tilde{L}_{n}\rightarrow \tilde{L}_{n}$
 by $$\tilde{D}(x+f)=Dx+D^{\ast}f, ~ \forall x \in
L_n, f\in L_n^{\ast}\eqno(3.7),$$ where $D^{\ast}f=-fD$. Then $\tilde{D}$  is an invertible, and by the direct computation, we have  $\tilde{D}\in Der_B(\tilde{L}_{n})$. Hence the metric $3$-Lie algebra $(\tilde{{L}}_{n},B)$ admits
a symplectic structure $\omega$ as follows
$$\omega(x+f, y+g)=B(\tilde{D}(x+f), y+g)=-f(Dy)+g(Dx).
\eqno(3.8)$$

\vspace{1mm} {\bf Remark 2 } By above discussion, from an arbitrary $3$-Lie algebra, we can construct
infinitely many metric sympletric $3$-Lie algebras.

\section{  Symplectic
structures  of $\mathrm{T}_{\theta}^{\ast}$-extensions}
\setcounter{equation}{0}
\renewcommand{\theequation}
{4.\arabic{equation}}

In papers \cite{bai4, bai5}, authors studied the extensions and module-extensions of
$3$-Lie algebras. In this  section we need
$T^*_{\theta}$-extension of $3$-Lie algebras to describe the
symplectic structures.

  \vspace{2mm}{\bf Lemma 4.1}
{\small\cite{bai4}} {\it Let $A$ be a $3$-Lie algebra over a
 field $F$, $A^{\ast}$ be the dual space of $A$, $\theta: A\wedge A\wedge A \rightarrow A^{\ast}$ be a linear mapping satisfying
$$\theta([x,u,v],y,z)+\theta([y,u,v],z,x)+\theta(x,y,[z,u,v])=\theta([x,y,z],u,v). \eqno(4.1 )$$
Then
$T_{\theta}^{\ast} A=A\oplus A^{\ast}$  is a $3$-Lie algebra in the following multiplication
$$
  [x+f, y+g, z+h]=[x, y, z]_{A}+\theta(x, y, z)+ ad^*(y, z)f+ ad^*(z, x)g+ad^*(x, y)h,\eqno(4.2)
$$
where $x, y, z\in A, f, g, h\in A^{\ast}$. The $3$-Lie algebra
$T_{\theta}^{\ast} A$ is called the $T_{\theta}^{\ast}$-extension of the $3$-Lie algebra $A$
by means of $\theta$.

If further, $\theta$ satisfies
$$\theta(x_{1}, x_{2}, x_{3})(x_{4})+\theta(x_{1}, x_{2}, x_{4})(x_{3})=0, \eqno(4.3)$$
for all $x_{1}, x_{2}, x_{3}, x_{4}\in A$, then the symmetric
bilinear form $B$ on $T_{\theta}^{\ast}A$ given
by
$$ B(x+f, y+g)=f(y)+g(x), ~ x, y\in A, f, g\in A^{\ast}, \eqno(4.4)$$
 defines
a metric structure on $T_{\theta}^{\ast}A$.}

\vspace{2mm}{\bf Theorem 4.2} { Let $A$ be a $3$-Lie algebra
admitting an invertible derivation  $D$, and $\theta: A\wedge A\wedge A \rightarrow A^{\ast}$ be a linear mapping satisfying Eqs.(4.1) and (4.3).
If there exists a linear mapping $\Psi: A\wedge A\rightarrow F$
satisfying for $x, y, z, u\in A$,
$$
  \Theta(x, y, z,u)=-(\Psi(x, [y ,z, u])-\Psi(y, [x, z, u])+\Psi(z, [x, y, u])-\Psi(u,[x, y, z])),\eqno(4.5)
 $$
 where $$\Theta(x, y, z,u)=\theta(Dx, y, z)u-\theta(Dy, z,u)x+\theta(Dz,u, x)y-\theta(Du, x, y)z,\eqno(4.6)$$
  then the metric $3$-Lie algebra $T_{\theta}^{\ast}A$ admits a symplectic structure.}

\vspace{1mm}{\bf Proof.} Let $B$ be the metric on the $3$-Lie algebra
$T_{\theta}^{\ast}A$ defined in Eq.(4.4). By Theorem
3.1, it suffices to prove that the existence of an invertible
skew-symmetric derivation of the metric $3$-Lie algebra $(T_{\theta}^{\ast}A, B)$.

Defines a linear mappings $H: A\rightarrow A^{\ast}$ and $\bar{D}: T_{\theta}^{\ast} A\rightarrow T_{\theta}^{\ast} A$, respectively,  by
$$B(Hx, y)= \Psi(x,y), ~ \forall x,y\in A,$$
 and $$\bar{D}(x+f)=Dx-Hx-f D, ~ \forall x\in A, f\in A^{\ast}.$$

  It is straightforward to see that $\bar{D}$ is invertible, since $D$ is so. And

\vspace{2mm}
 $B(\bar{D}(x+f), y+g)=B(Dx-Hx-f D, y+g)=g(Dx)-f(Dy)-F(x,y),$

 \vspace{2mm} $B(x+f, \bar{D}( y+g))=B(x+f, Dy-Hy-g D)=-g(Dx)+f(Dy)-F(y,x).$

 Therefore, $\bar{D}$ is  skew-symmetric with respect to the metric  $B$.

Further, since $D$ is a derivation of $A$, for $x, y, z\in A$ and $f, g, h\in A^*$  we get

\vspace{2mm}
$[\bar{D}(x+f),y+g,z+h]+[x+f,\bar{D}(y+g),z+h]$

\vspace{2mm}$+[x+f,y+g,\bar{D}(z+h)]-\bar{D}[x+f,y+g,z+h]$

\vspace{2mm}\noindent
$=[Dx-Hx-fD,y+g,z+h]+[x+f,Dy-Hy-gD,z+h]$

\vspace{2mm}$+[x+f,y+g,Dz-Hz-hD]-\bar{D}([x,y,z]+\theta(x,y,z)$

\vspace{2mm}$+ad^*(y,z)f+ad^*(z,x)g+ad^*(x,y)h)$

\vspace{2mm}\noindent$=[Dx,y,z]+\theta(Dx,y,z)-ad^*(y,z)(Hx+fD)+ad^*(z, Dx)g$

\vspace{2mm}$+ad^*(Dx, y)h+[x, Dy ,z]+\theta(x, Dy, z)+ad^*(Dy, z)f$

\vspace{2mm}$-ad^*(z, x)(Hy+gD)+ad^*(x, Dy)h+[x, y, Dz]+\theta(x, y, Dz)$

\vspace{2mm}$+ad^*(y, Dz)f+ad^*(Dz, x)g-ad^*(x, y)(Hz+hD)-D[x, y, z]$

\vspace{2mm}$+H[x, y, z]+\theta(x, y, z)D-D^*ad^*(y, z)f-D^*ad^*(z, x)g-D^*ad^*(x, y)h $

\vspace{2mm}\noindent$=\theta(Dx, y, z)+\theta(x, Dy, z)+\theta(x, y, Dz)+\theta(x, y, z)D-ad^*(y ,z)Hx$

\vspace{2mm}$-ad^*(z, x)Hy-ad^*(x, y)Hz+H[x, y, z].$

From  Eqs.(4.5) and (4.6)  and $\Psi(x,y)=B(Hx,y)=Hx(y)$ for all $x,y\in A$,  for arbitrary $
u\in A$,

\vspace{2mm}
 $\theta(Dx,y,z)u+\theta(x,Dy,z)u+\theta(x,y,Dz)u+\theta(x,y,z)
Du$

\vspace{2mm}$+B(Hx,[y,z,u])+B(Hy,[z,x,u])+B(Hz,[x,y,u])+B(H[x,y,z],u)$

\vspace{2mm}\noindent$=\Theta(x,y,z,u)+\Psi(x, [y ,z, u])-\Psi(y, [x, z, u])+\Psi(z, [x, y,u])-\Psi(u,[x, y, z]=0.$

Therefore, $\bar{D}$ is an invertible derivation of $T_{\theta}^{\ast}A$. The proof is completed. $\Box$

\vspace{2mm}{\bf Lemma 4.3} {\it Let $A $ be a nilpotent
$3$-Lie algebra over  $F$, $I$ be a nonzero ideal of $A$. Then
$I\cap Z(A)\neq 0$.}

\vspace{2mm}{\bf Proof}. If $A$ is abelian, the result is evident.

If $A$ is non-abelian, and  $I$ is a nonzero ideal of $A$. Then for every $x, y\in A$, the left multiplication $ad(x, y): A\rightarrow A$ is nilpotent (\cite{K1}). Therefore,
 the inner derivation algebra $ad (A)$ of the $3$-Lie algebra $A$
is constituted by nilpotent mappings. Since $ad(x, y)(I)\subseteq I$, for all $x, y\in A$, by Theorem 3.3 in \cite{Hym}, there exists non-zero element $z\in I$ such that $ad(x, y)(z)=0, \forall x, y\in L$. Therefore, $z\in I\cap Z(A)$. $\Box$

\vspace{2mm}{\bf Lemma 4.4} {\it Let $(A, B)$ be a non-abelian nilpotent metric
$3$-Lie algebra over $F$. Then there exists a non-zero isotropic ideal of $A$.}

\vspace{2mm}{\bf  Proof. } Denotes $J= A^1\cap Z(A)$. By Lemma 4.3 $J$ is a non-zero ideal of $A$. Thanks to Lemma 2.3 in paper \cite{bai3}, $Z(A)^{\bot}=A^1=[A, A, A].$ Then, $J\subseteq J^{\bot},$ that is, $J$ is a non-zero isotropic
ideal of $A$. $\Box$

\vspace{2mm}{\bf Lemma 4.5}{\cite{bai4} {\it  Let $(L, B)$ be a nilpotent metric $3$-Lie algebra of dimension $m$. If
$J$ is an isotropic ideal of $L$, then $L$ contains a maximally isotropic ideal I of dimension
$[\frac{m}{2}]$ containing $J$. Moreover,

1) If $m$ is even, then $L$ is isometric to some $T^*_{\theta}$-extension
of $L/I$.

2) If $m$ is odd, then the ideal $I^{\bot}$
is an abelian ideal of $L$, and $L$ is isometric to a
non-degenerate ideal of codimension one in some $T^*_{\theta}$-extension of $L/I$.}

\vspace{2mm}{\bf Theorem 4.6} {\it Let $(L, B)$ be a  non-abelian  nilpotent metric
3-Lie algebra over an algebraically closed field $F$
which admits a skew-symmetric invertible derivation $\bar{D}$.
Then there exists a $3$-Lie algebra $A$, an invertible derivation
$D$ of $A$ and  $\theta: A\wedge A\wedge A\rightarrow A^*$
satisfying Eq.(4.1) such that $L=T_{\theta}^{\ast}A$. And There exists $\Psi: A\wedge A \rightarrow F$  such that
$\Theta(x,y,z,u)$ defined by Eq.(4.6) satisfying Eq.(4.5).}

\vspace{1mm}{\bf Proof.} By Lemma 4.3 and Lemma 4.4, $I=L^1\cap Z(L)$
is a non-zero isotropic characteristically ideal of the $3$-Lie algebra $L$.  From Theorem 3.1,
there exists a non-degenerate skew-symmetric bilinear form $\omega$ on $L$ such that the invertible derivation
$\bar{D}$ satisfies $$\omega(\bar{D}x, y)+\omega(x, \bar{D}y)=0.$$  Therefore, the dimension of the $3$-Lie algebra $L$ is even.

 Since the $3$-Lie algebra $L$ is nilpotnet, the inner derivation algebra $Ad(L)$
is a nilpotent Lie algebra. Then the Lie algebra $T=Ad(L)\oplus F\bar{D}$ is solvable.  By Lemma 3.2 in \cite{Bor} and Lemma 4.5,
 there exists a maximal
isotropic ideal $J$ containing the isotropic ideal $I=L^1\cap Z(L)$,
 and $\theta: (L/J)\wedge (L/J)\wedge (L/J)\rightarrow (L/J)$
satisfying Eq.(4.1)  such that  the metric $3$-Lie algebra $(L, B)$ is
isomorphic to the $T_{\theta}^{\ast}$-extension $T^*_{\theta}(L/J)$, and  $J$ is stable by $\bar{D}.$
Let $J'$ be a complement of $J$ in the vector space $L$, that is, $L=J'\oplus J$.
Then for every $x\in J, ~ y\in J'$, we have  $\bar{D}(x)\in J$ and $\bar{D}(y)=y_1+y_2$, where $y_1\in J'$ and $y_2\in J$.
 Denotes the $3$-Lie algebra $L/J$ by $A$. Then $A^*$ is isomorphic to $J$ as subspaces and it is stable by $\bar{D}$.

 Therefore, we can define linear
mappings $D_{11}: A\rightarrow A,$ ~ $D_{21}: A\rightarrow A^{\ast}$, and
$D_{22}: A^{\ast}\rightarrow A^{\ast}$ by
$$ \bar{D}(x+f)=D_{11}x+D_{21}x+D_{22}f,~~  \forall x\in
A, f\in A^{\ast}.\eqno(4.7)$$

Clearly, $D_{11}$ and $D_{22}$ must be
invertible since $\bar{D}$ is so. And for every $x, y\in A, f, g\in A^*$
 \begin{eqnarray}
 0&=&B(\bar{D}(x+f),y+g)+B(x+f,\bar{D}(y+g))\nonumber\\
 &=&B(D_{11}x+D_{21}x+D_{22}f,y+g)+B(x+f,D_{11}y+D_{21}y+D_{22}g)\nonumber\\
 &=&g(D_{11}x)+D_{21}x(y)+D_{22}f(y)+f(D_{11}y)+D_{21}y(x)+D_{22}g(x).\nonumber\hspace{3cm}\hfill(4.8)
 \end{eqnarray}

From the above equation, we obtain that in the case $x=0,$ $ g=0,$
$$ D_{22}f(y)=-fD_{11}(y), ~~ \forall
y\in A, ~f\in a^{\ast}, $$
and in the case $f=g=0$,
$$ B(D_{21}x, y)+B(D_{21}y, x)=0, \forall x, y\in A.$$

Let $H=-D_{21}: A\rightarrow A^*$ and $D=D_{11}: A\rightarrow A$. Since $\bar{D}$ is a derivation of $L$, by Eq.(4.2)

\vspace{2mm}\hspace{4mm}$
0=[\bar{D}x,y,z]+[x,\bar{D}y,z]+[x,y,\bar{D}z]-\bar{D}[x,y,z]$

\vspace{2mm}$=[Dx-Hx,y,z]+[x,Dy-Hy,z]+[x,y,Dz-Hz]-\bar{D}([x,y,z]+\theta(x,y,z))$

\vspace{2mm}$=[Dx,y,z]+\theta(Dx,y,z)-ad^*(y, z)Hx+[x,Dy,z]+\theta(x,Dy,z)-ad^*(z, x)Hy$

\vspace{2mm}\hspace{4mm}$+[x,y,Dz]+\theta(x,y,Dz)-ad^*(x, y)Hz-D[x,y,z]+H[x,y,z]+\theta(x,y,z)D$

\vspace{2mm}$=[Dx,y,z]+[x,Dy,z]+[x,y,Dz]-D[x,y,z]$

\vspace{2mm}\hspace{4mm}$+\theta(Dx,y,z)+\theta(x,Dy,z)+\theta(x,y,Dz)+\theta(x,y,z)D $

\vspace{2mm}\hspace{4mm}$-ad^*(y, z)Hx-ad^*(z, x)Hy-ad^*(x, y)Hz+H[x,y,z], ~\forall x,y,z\in A.$

Therefore, we have
$$[Dx,y,z]+[x,Dy,z]+[x,y,Dz]-D[x,y,z]=0,~~  ~\forall x,y,z\in A, \eqno(4.9)
$$

\hspace{2.5cm}$\theta(Dx,y,z)+\theta(x,Dy,z)+\theta(x,y,Dz)+\theta(x,y,z)D$

\vspace{2mm}\hspace{1.8cm}$= ad^*(y, z)Hx+ad^*(z, x)Hy+ad^*(x, y)Hz-H[x,y,z], ~\forall x,y,z\in A.$
\hfill(4.10)

Therefore, $D$ is an invertible derivation of $A$. Denotes
\begin{equation*}
   \Theta(x,y,z,u)=\theta(Dx,y,z)u-\theta(Dy,z,u)x+\theta(Dz,u,x)y-\theta(Du,x,y)z,~~\forall x, y, z, u\in A.
\end{equation*}

Defines bilinear mapping $\Psi: A\times A \rightarrow F$ by
$$\Psi(x,y)=-B(Hx,y)=-Hx(y), ~~ \forall x, y\in A.$$
Then $\Psi$ is skew-symmetric and satisfies $ ~\forall x, y, z,\omega\in A,$
$$
  \Theta(x,y,z,\omega)+(\Psi(x,[y,z,\omega])-\Psi(y,[x,z,\omega])+\Psi(z,[x,y,\omega])-\Psi(\omega,[x,y,z]))=0.$$ The result follows. $\Box$

\vspace{2mm}The following result gives a characterization of 3-Lie
algebras admitting an invertible derivation. Note that the result
is valid for an arbitrary base field of characteristic zero (not
necessarily algebraically closed).

\vspace{2mm}{\bf Theorem 4.7 } {\it Let  $A$ be a $3$-Lie algebra over
a field $F$ with a characteristic zero. Then there exists an invertible derivation $D$ of $A$ if and only if $A$ is
isomorphic to the quotient $3$-Lie algebra $L/J$ of a
metric symplectic $3$-Lie algebra $(L, B, \omega)$ by a lagrangian and
completely isotropic ideal $J$.}

\vspace{1mm}{\bf Proof.} If $A$ admits an invertible derivation. From Theorem 4.2, let $\theta=0, \Psi=0, H=0$
then the $3$-Lie algebra $L=A\oplus A^*$  obtained by
$T^{\ast}_{0}$-extension of $A$ is a metric symplectic $3$-Lie
algebra.

We define
$$\bar{D}:L\rightarrow
L, ~~ \bar{D}(x+f)=Dx-fD, ~ \forall x\in
A, f\in A^{\ast},$$ and
$$\omega(x+f, y+g)=B(\bar{D}(x+f), y+g)=g(Dx)-f(Dy), ~ x, y\in A, f, g\in A^{\ast}.$$
 Then $J=A^{\ast}$ is a
lagrangian ideal of the symplectic
$3$-Lie algebra $(L, \omega)$, and is a completely isotropic ideal of the metric
$3$-Lie algebra $(L, B)$, and $A$ is isomorphic to the quotient $3$-Lie algebra $L/J$.

Conversely, suppose that the $3$-Lie algebra $A$ is isomorphic to $L/J$,
where $(L, B, \omega)$ is a metric symplectic $3$-Lie algebra and
$J$ is a lagrangian completely isotropic ideal of $L$.
By Theorem 3.4 in \cite{bai4}, $L$ is isometrically isomorphic
to $T_{\theta}^{\ast}(L/J)=T_{\theta}^{\ast}A$ since
$J$ is completely isotropic. From Theorem 3.1, there exists  a skew-symetric invertible derivation $\bar{D}$ of the metric $3$-Lie algebra $(L, B).$
From Eq.(3.2), $\bar{D}(J)=J.$ Then by
 the same argument used
in the proof of Theorem 4.6, the projection
$\bar{D}|_ A: A\rightarrow A$ provides a non-singular derivation of
$A$. $\Box$

\vspace{3mm} At last of the paper, we give the characterization of metric symplectic double extensions of $3$-Lie algebras.

\vspace{2mm}{\bf Lemma 4.8}\cite{bai5} {\it Let $(A, B)$ be a metric $3$-Lie algebra,
$b$ be another $3$-Lie algebra and $\pi=ad^{\ast}: b\times
b\rightarrow End(b^{\ast})$ be the coadjoint representation of
$b$. Suppose that $(A, \psi)$ is a representation of $b$, where $\psi: b\wedge
b\rightarrow  End(A)$ satisfies $\psi(b, b)\subseteq Der_{B}(A)$.
Let $\tilde{A}=b^{\ast}\oplus A\oplus b,$ ~and  $\phi: A\otimes
A\otimes b\rightarrow b^{\ast}$ defined by  for any $x_{1}, x_{2}\in A, y, z\in
b$
$$
  \phi(x_{1},x_{2},y)(z)=-\phi(x_{2},x_{1},y)(z)=B(\psi(y,z)x_{1},x_{2}).
$$
If $\psi$ satisfies $\psi(b^{1}, b)(A)=\psi(b,b)(A^{1})=0.$
Then $(\tilde{A}, T)$ is a metric $3$-Lie algebra in the following multiplication,
 $\forall y_{1}, y_{2}, y_{3}\in b,$ $\forall  x_{1}, x_{2}, x_{3}\in A,$
$ \forall f_{1}, f_{2}, f_{3}\in b^{\ast}$,

\vspace{2mm}\hspace{6mm}$[y_{1}+x_{1}+f_{1}, y_{2}+x_{2}+f_{2}, y_{3}+x_{3}+f_{3}]$

  \vspace{2mm}$=[y_{1}, y_{2}, y_{3}]_{b}+[x_{1}, x_{2}, x_{3}]_{A}+\psi(y_{2}, y_{3})x_{1}-\psi(y_{1}, y_{3})x_{2}+\psi(y_{1}, y_{2})x_{3}
  +\pi(y_{2}, y_{3})f_{1}$

  \vspace{2mm}\hspace{4mm}$ -\pi(y_{1}, y_{3})f_{2}+\pi(y_{1}, y_{2})f_{3}+\phi(x_{1}, x_{2}, y_{3})-\phi(x_{1}, x_{3}, y_{2})+\phi(x_{2}, x_{3}, y_{1}).$\hfill(4.11)

$$T(y_{1}+x_{1}+f_{1},y_{2}+x_{2}+f_{2})=B(x_{1},x_{2})+f_{1}(y_{2})+f_{2}(y_{1}). \eqno(4.12)$$
$\Box$}

\vspace{1mm}  In Lemma 4.8,  if $b=Fe_{1}+Fe_{2}$ is a two-dimensional $3$-Lie
algebra, then  $$\psi: b\wedge b\rightarrow A, ~ \psi(e_{1},e_{2})=\delta\in Der_B(A).$$
Therefore, $\phi: A\otimes
A\otimes b\rightarrow b^{\ast}$ defined by  for any $x_{1}, x_{2}\in A, e_1, e_2\in
b$
$$
  \phi(x_{1}, x_{2}, e_1)(e_2)=-\phi(x_{2}, x_{1}, e_1)(e_2)=B(\psi(e_1, e_2)x_{1}, x_{2})=B(\delta x_1, x_2),\eqno(4.13)
$$
$$
  \phi(x_{1}, x_{2}, e_2)(e_1)=-B(\delta x_1, x_2), ~~  \phi(x_{1},x_{2},e_1)(e_1)=\phi(x_{1},x_{2},e_2)(e_2)=0.\eqno(4.13')
$$
Then we say  that $(\tilde{A}=Fe_{1}+Fe_{2}\oplus A\oplus
Fe_{1}^{\ast}+Fe_{2}^{\ast},~~  T)$ is
the double extension of $A$ by means of the derivation
$\psi(e_{1},e_{2})=\delta$, and the multiplication is for $\forall x, y, x\in
A,$ $\alpha, \alpha^{'}, $ $\beta, \beta^{'},$ $\gamma_{1}, \gamma_{1}^{'},\gamma_{2},$ $\gamma_{2}^{'}\in F,$ $ e_1^*, e_2^*\in b^* $ ( where  $e_i^*(e_j)=\delta_{ij}, $ $1\leq 1, j\leq 2 $),

\vspace{2mm}$[\alpha e_{1}+x+\alpha^{'}e_{1}^{\ast},\beta
e_{2}+y+\beta^{'}e_{2}^{\ast},\gamma_{1}e_{1}+\gamma_{2}e_{2}+z+\gamma_{1}^{'}e_{1}^{\ast}+\gamma_{2}^{'}e_{2}^{\ast}]$

\vspace{2mm}\noindent$=[x,y,z]+\delta(-\beta\gamma_{1}
x-\alpha\gamma_{2}y+\alpha\beta
z)+\phi(x,y,\gamma_{1}e_{1}+\gamma_{2}e_{2})-\phi(x,z,\beta
e_{2})+\phi(y,z,\alpha e_{1}),$\hfill(4.14)

 \vspace{2mm}\noindent and the metric is $$T(\alpha
e_{1}+x+\alpha^{'}e_{1}^{\ast},\beta
e_{2}+y+\beta^{'}e_{2}^{\ast})=B(x,y)+\alpha\beta^{'}+\beta\alpha^{'}.\eqno(4.15)$$

By the above notations we have the following result.

\vspace{2mm}{\bf Theorem 4.9} {\it Let $(A, B)$ be a metric $3$-Lie algebra, $D$ be an invertible derivation of $A$ and $D\in Der_B(A)$, and $\delta\in Der_B(A)$ satisfy
$$\delta D-D\delta=2\delta.\eqno(4.16)$$
Let $(\tilde{A}=b\oplus A\oplus b^*, T)$  be the double extension  of $A$ by means of the derivation
$\delta$, where $b=Fe_1+Fe_2$ be the $2$-dimensional $3$-Lie algebra. Then the linear endomorphism $\tilde{D}$ of $\tilde{A}$ defined by
$$
\tilde{D}|_{A}=D, ~\tilde{D}{e_{i}}=-e_{i},~ \tilde{D}{e_{i}^{\ast}}=e_{i}^{\ast},
~i=1, 2 \eqno(4.17)
$$
is an invertible derivation of the $3$-Lie algebra $(\tilde{A}, T)$, and $\tilde{D}\in Der_T(\tilde{A})$.}

{\bf Proof } Let $\psi: b\wedge b\rightarrow A, ~ \psi(e_{1},e_{2})=\delta\in Der_B(A).$ By the above discussion, $(\tilde{A}=b\oplus A\oplus b^*, T)$ is the double extension  of $A$ by means of the derivation
$\delta$.

By Eq.(4.17), the linear mapping $\tilde{D}: \tilde{A} \rightarrow \tilde{A} $ is invertible.  From Lemma 4.8 and Eq.(4.14), $\forall x, y, z\in
A, \alpha, \alpha^{'}, \beta, \beta^{'},\gamma_{1}, \gamma_{1}^{'},\gamma_{2},\gamma_{2}^{'}\in F,$

\vspace{2mm}\hspace{6mm}$\tilde{D}[\alpha e_{1}+x+\alpha^{'}e_{1}^{\ast},\beta
e_{2}+y+\beta^{'}e_{2}^{\ast},\gamma_{1}e_{1}+\gamma_{2}e_{2}+z+\gamma_{1}^{'}e_{1}^{\ast}+\gamma_{2}^{'}e_{2}^{\ast}]$

$=D[x,y,z]+D\delta(-\beta\gamma_{1}
x-\alpha\gamma_{2}y+\alpha\beta
z)+\phi(x,y,\gamma_{1}e_{1}+\gamma_{2}e_{2})$

\vspace{2mm}\hspace{6mm}$-\phi(x,z,\beta
e_{2})+\phi(y,z,\alpha e_{1}).$

\vspace{2mm}Thanks to Eqs.(4.16) and (4.17),

\vspace{2mm}\hspace{6mm}$[\tilde{D}(\alpha e_{1}+x+\alpha^{'}e_{1}^{\ast}),\beta
e_{2}+y+\beta^{'}e_{2}^{\ast},\gamma_{1}e_{1}+\gamma_{2}e_{2}+z+\gamma_{1}^{'}e_{1}^{\ast}+\gamma_{2}^{'}e_{2}^{\ast}]$

\vspace{2mm}\hspace{6mm}$+[\alpha e_{1}+x+\alpha^{'}e_{1}^{\ast},\tilde{D}(\beta
e_{2}+y+\beta^{'}e_{2}^{\ast}),\gamma_{1}e_{1}+\gamma_{2}e_{2}+z+\gamma_{1}^{'}e_{1}^{\ast}+\gamma_{2}^{'}e_{2}^{\ast}]$

\vspace{2mm}\hspace{6mm}}$+[\alpha e_{1}+x+\alpha^{'}e_{1}^{\ast},\beta
e_{2}+y+\beta^{'}e_{2}^{\ast},\tilde{D}(\gamma_{1}e_{1}+\gamma_{2}e_{2}+z+\gamma_{1}^{'}e_{1}^{\ast}+\gamma_{2}^{'}e_{2}^{\ast})]$

\vspace{2mm}$=[-\alpha e_{1}+Dx+\alpha^{'}e_{1}^{\ast},\beta
e_{2}+y+\beta^{'}e_{2}^{\ast},\gamma_{1}e_{1}+\gamma_{2}e_{2}+z+\gamma_{1}^{'}e_{1}^{\ast}+\gamma_{2}^{'}e_{2}^{\ast}]$

\vspace{2mm}\hspace{6mm}$+[\alpha e_{1}+x+\alpha^{'}e_{1}^{\ast},-\beta
e_{2}+Dy+\beta^{'}e_{2}^{\ast},\gamma_{1}e_{1}+\gamma_{2}e_{2}+z+\gamma_{1}^{'}e_{1}^{\ast}+\gamma_{2}^{'}e_{2}^{\ast}]$

\vspace{2mm}\hspace{6mm}$+[\alpha e_{1}+x+\alpha^{'}e_{1}^{\ast},\beta
e_{2}+y+\beta^{'}e_{2}^{\ast},-\gamma_{1}e_{1}-\gamma_{2}e_{2}+Dz+\gamma_{1}^{'}e_{1}^{\ast}+\gamma_{2}^{'}e_{2}^{\ast}]$

\vspace{2mm}$=[Dx,y,z]+[x,Dy,z]+[x,y,Dz]$

\vspace{2mm}\hspace{6mm}$+\delta D(-\beta\gamma_{1} x-\alpha\gamma_{2}y+\alpha\beta z)-2\delta(-\beta\gamma_{1} x-\alpha\gamma_{2}y+\alpha\beta z)$

\vspace{2mm}\hspace{6mm}$+\phi(Dx, y,\gamma_{1}e_{1}+\gamma_{2}e_{2})
-\phi(Dx,z,\beta
e_{2})+\phi(y,z,-\alpha e_{1})$

\vspace{2mm}\hspace{6mm}$+\phi(x,Dy,\gamma_{1}e_{1}+\gamma_{2}e_{2})-\phi(x,z,-\beta
e_{2})+\phi(Dy,z,\alpha e_{1})$

\vspace{2mm}\hspace{6mm}$+\phi(x,y,-\gamma_{1}e_{1}-\gamma_{2}e_{2})-\phi(x,Dz,\beta
e_{2})+\phi(y,Dz,\alpha e_{1})$

\vspace{2mm}$=D[x,y,z]+D\delta (-\beta\gamma_{1} x-\alpha\gamma_{2}y+\alpha\beta z)$

\vspace{2mm}\hspace{6mm}$+\phi(Dx, y,\gamma_{1}e_{1}+\gamma_{2}e_{2})
-\phi(Dx,z,\beta
e_{2})-\phi(y,z,\alpha e_{1})$

\vspace{2mm}\hspace{6mm}$+\phi(x,Dy,\gamma_{1}e_{1}+\gamma_{2}e_{2})+\phi(x,z,\beta
e_{2})+\phi(Dy,z,\alpha e_{1})$

\vspace{2mm}\hspace{6mm}$-\phi(x,y,\gamma_{1}e_{1}+\gamma_{2}e_{2})-\phi(x,Dz,\beta
e_{2})+\phi(y,Dz,\alpha e_{1}). $

 From Eqs.(4.13) and (4.16),

\vspace{2mm}\hspace{6mm}$(\phi(Dx, y,\gamma_{1}e_{1}+\gamma_{2}e_{2})+\phi(x,Dy,\gamma_{1}e_{1}+\gamma_{2}e_{2})-\phi(x,y,\gamma_{1}e_{1}+\gamma_{2}e_{2}))(e_{1})
$

\vspace{2mm}$=B(-\gamma_{2}\delta Dx,y)+B(-\gamma_{2}\delta x,Dy)+B(\gamma_{2}\delta x,y)
$

\vspace{2mm}$=B(-\gamma_{2}\delta Dx,y)+B(\gamma_{2}D\delta x, y)+B(\gamma_{2}\delta x,y)
$

\vspace{2mm}$=-\gamma_2B((\delta D-D\delta-2\delta)x, y)-B(\gamma_{2}\delta x,y)$

\vspace{2mm}$=B(-\gamma_{2}\delta x,y)=\phi(x,y,\gamma_{1}e_{1}+\gamma_{2}e_{2}))(e_{1}),$

\vspace{2mm}\hspace{6mm}$(\phi(Dx, y,\gamma_{1}e_{1}+\gamma_{2}e_{2})+\phi(x,Dy,\gamma_{1}e_{1}+\gamma_{2}e_{2})-\phi(x,y,\gamma_{1}e_{1}+\gamma_{2}e_{2}))(e_{2})
$

\vspace{2mm}$=B(\gamma_{1}\delta Dx,y)+B(\gamma_{1}\delta x,Dy)-B(\gamma_{1}\delta x,y)
$

\vspace{2mm}$=\gamma_{1}B((D\delta-D\delta-2\delta)x,y)+B(\gamma_{1}\delta x, y)$

\vspace{2mm}$=B(\gamma_{1}\delta x,y)=\phi(x,y,\gamma_{1}e_{1}+\gamma_{2}e_{2}))(e_{2})$.

Then we have

$\phi(Dx, y,\gamma_{1}e_{1}+\gamma_{2}e_{2})+\phi(x,Dy,\gamma_{1}e_{1}+\gamma_{2}e_{2})-\phi(x,y,\gamma_{1}e_{1}+\gamma_{2}e_{2})=\phi(y,z,\gamma_{1}e_{1}+\gamma_{2}e_{2}).$

\vspace{2mm}Similarly,

\vspace{2mm}$-\phi(Dx,z,\beta e_{2})+\phi(x,z,\beta e_{2})-\phi(x,Dz,\beta e_{2})=-\phi(x,z,\beta e_{2}),$

\vspace{2mm}$-\phi(y,z,\alpha e_{1})+\phi(Dy,z,\alpha e_{1})+\phi(y,Dz,\alpha e_{1})=\phi(y,z,\alpha e_{1}).$

Therefore, $\tilde{D}$ satisfies

 $\tilde{D}[\alpha e_{1}+x+\alpha^{'}e_{1}^{\ast},\beta
e_{2}+y+\beta^{'}e_{2}^{\ast},\gamma_{1}e_{1}+\gamma_{2}e_{2}+z+\gamma_{1}^{'}e_{1}^{\ast}+\gamma_{2}^{'}e_{2}^{\ast}]$

\vspace{2mm}$=[\tilde{D}(\alpha e_{1}+x+\alpha^{'}e_{1}^{\ast}),\beta
e_{2}+y+\beta^{'}e_{2}^{\ast},\gamma_{1}e_{1}+\gamma_{2}e_{2}+z+\gamma_{1}^{'}e_{1}^{\ast}+\gamma_{2}^{'}e_{2}^{\ast}]$

\vspace{2mm}$+[\alpha e_{1}+x+\alpha^{'}e_{1}^{\ast},\tilde{D}(\beta
e_{2}+y+\beta^{'}e_{2}^{\ast}),\gamma_{1}e_{1}+\gamma_{2}e_{2}+z+\gamma_{1}^{'}e_{1}^{\ast}+\gamma_{2}^{'}e_{2}^{\ast}]$

\vspace{2mm}$+[\alpha e_{1}+x+\alpha^{'}e_{1}^{\ast},\beta
e_{2}+y+\beta^{'}e_{2}^{\ast},\tilde{D}(\gamma_{1}e_{1}+\gamma_{2}e_{2}+z+\gamma_{1}^{'}e_{1}^{\ast}+\gamma_{2}^{'}e_{2}^{\ast})].$

Again  by Eqs.(4.15) and (4.17),

\vspace{2mm}\hspace{6mm}$T(\tilde{D}(\alpha e_{1}+\beta e_2+\epsilon x+\alpha^{'}e_{1}^{\ast}+\beta'e_2^*), \lambda e_1+\mu e_2+\nu y+
\lambda^{'}e_{1}^{\ast}+\nu^{'}e_{2}^{\ast})$

\vspace{2mm}\hspace{6mm}$+T(\alpha e_{1}+\beta e_2+\epsilon x+\alpha^{'}e_{1}^{\ast}+\beta'e_2^*, \tilde{D}( \lambda e_1+\mu e_2+\nu y+
\lambda^{'}e_{1}^{\ast}+\nu^{'}e_{2}^{\ast}))$

\vspace{2mm}$=T(-\alpha e_{1}-\beta e_2+\epsilon Dx+\alpha^{'}e_{1}^{\ast}+\beta'e_2^*,  \lambda e_1+\mu e_2+\nu y+
\lambda^{'}e_{1}^{\ast}+\nu^{'}e_{2}^{\ast})$

\vspace{2mm}\hspace{6mm}$+T(\alpha e_{1}+\beta e_2+\epsilon x+\alpha^{'}e_{1}^{\ast}+\beta'e_2^*, -\lambda e_1-\mu e_2+\nu Dy+
\lambda^{'}e_{1}^{\ast}+\nu^{'}e_{2}^{\ast})$

\vspace{2mm}$=B(\epsilon Dx, \nu y)+B(\epsilon x, \nu Dy)-\alpha \lambda'-\beta \mu'+ \alpha'\lambda+\beta'\mu+\alpha \lambda'+\beta \mu'-\alpha'\lambda-\beta'\mu=0.$

Summarizing above discussion, we obtain that
$\tilde{D}$ is an invertible  derivation of the metric $3$-Lie algebra $(\tilde{A}, T)$ and $\tilde{D}\in Der_T(\tilde{A})$. $\Box$

If $(A, B)$ be a metric $3$-Lie algebra and $D\in Der_B(A)$ is invertible. From  Eq.(3.2), $(A, B, \omega)$ is a metric symplectic $3$-Lie
algebra, where $\omega(x, y)=B(Dx, y), \forall x, y\in A$. Then  we obtain the following result.

\vspace{2mm}{\bf Corollary} {\it Let $(A, B)$ be a metric $3$-Lie algebra, $D$ be an invertible derivation of $A$, $D\in Der_B(A)$ and
$\delta\in Der_B(A)$ satisfy Eq.(4.16). Then
the $3$-Lie algebra $(\tilde{A}, T, \tilde{\omega})$ is a metric symplectic $3$-Lie
algebra,  which is called the metric symplectic double extension of
$(A, B,\omega)$, where $(\tilde{A}, T)$  is the double extension of $(A, B)$ by means of $\delta$, and $\tilde{\omega}$ is defined by
$$\tilde{\omega}(x, y)=\omega(x, y), ~ \tilde{\omega}(e_{1},e_{2}^{\ast})=\tilde{\omega}(e_{2}, e_{1}^{\ast})=-1, ~~ \forall  x, y \in A.\eqno(4.18)$$}

{\bf Proof } The result follows from Theorem 4.9 and Theorem 3.1, directly. $\Box$

\section*{Acknowledgements}
 The first author would like to thank
supports from National Natural Science Foundation of China
(11371245) and Natural Science Foundation of Hebei Province
(A2014201006), China.

\bibliography{}

\end{document}